\documentclass[twoside,11pt]{amsart}%{article}%
\usepackage{amsmath,amssymb,amsthm}

\usepackage{hyperref}

\begin{document}

\setlength{\textwidth}{145mm} \setlength{\textheight}{203mm}

\frenchspacing

%\pagestyle{myheadings}

%\markboth{Dimitar Mekerov}{Some conditions the curvature tensor of
%the Bismut connection ... }

% THEOREM Environments ---------------------------------------------------
\numberwithin{equation}{section}
\newtheorem{thm}{Theorem}[section]
\newtheorem{lem}[thm]{Lemma}
\newtheorem{prop}[thm]{Proposition}
\newtheorem{cor}[thm]{Corollary}
\newtheorem{probl}[thm]{Problem}

\newtheorem{defn}{Definition}[section]
\newtheorem{rem}{Remark}[section]
\newtheorem{exa}{Example}%[section]
\newtheorem{case}{Case}
%\theoremstyle{remark}

% MATH -------------------------------------------------------------------

\newcommand{\X}{\mathfrak{X}}
\newcommand{\s}{\mathfrak{S}}
\newcommand{\W}{\mathcal{W}}
\newcommand{\R}{\mathbb{R}}
\newcommand{\dd}{\mathrm{d}}
\newcommand{\n}{\nabla}
\newcommand{\nn}{\nabla'}
\newcommand{\lm}{\lambda}
\newcommand{\ta}{\theta}
\newcommand{\pd}{\partial}
\newcommand{\LL}{\mathcal{L}}

\newcommand{\de}{\mathrm{d}}
\newcommand{\im}{\mathrm{Im}}
\newcommand{\id}{\mathrm{id}}

\newcommand{\ie}{i.e. }
\newcommand{\norm}[1]{\left\Vert#1\right\Vert ^2}
\newcommand{\nnorm}[1]{\left\Vert#1\right\Vert ^{*2}}
\newcommand{\tr}{{\rm tr}}

\newcommand{\thmref}[1]{Theorem~\ref{#1}}
\newcommand{\propref}[1]{Proposition~\ref{#1}}
\newcommand{\corref}[1]{Corollary~\ref{#1}}
\newcommand{\secref}[1]{\S\ref{#1}}
\newcommand{\lemref}[1]{Lemma~\ref{#1}}
\newcommand{\dfnref}[1]{Definition~\ref{#1}}

%%% ----------------------------------------------------------------------

\title{Natural connections on conformal Riemannian $P$-manifolds}

\author{Dobrinka Gribacheva, Dimitar Mekerov}

%%% ----------------------------------------------------------------------
\maketitle
%%% ----------------------------------------------------------------------

{\small
%\begin{abstract}
\textbf{Abstract}

The class $\W_1$ of conformal Riemannian $P$-manifolds is the
largest class of Riemannian almost product manifolds, which is
closed with respect to the group of the conformal transformations
of the Riemannian metric. This class is an analogue of the class
of conformal K\"ahler manifolds in almost Hermitian geometry. In
the present work\footnote{Partially supported by project
NI11-FMI-004 of the Scientific Research Fund, Paisii Hilendarski
University of Plovdiv, Bulgaria} we study the natural connections
on the manifolds $(M, P, g)$ from the class $\W_1$, i.e. the
linear connections preserving the almost product structure $P$ and
the Riemannian metric $g$. We find necessary and sufficient
conditions the curvature tensor of such a connection to have
similar properties like the ones of the K\"ahler tensor in
Hermitian geometry. We determine the type of the manifolds
admitting a natural connection with a parallel torsion.


}
%\end{abstract}

%%%%%%%%%%%%%%%%%%%%%%%%%%%%%%%%%%%%%%%%%%%%%%%%%%%%%%%%%%%%%%%%%%%%%%%%%%%0

\section{Introduction}
K. Yano initiated in \cite{Ya} the study of Riemannian almost
product manifolds. A Riemannian almost product manifold $(M, P,
g)$ is a differentiable manifold $M$ with an almost product
structure $P$ and a Riemannian metric $g$ such that $P^2x = x$ and
$g(Px, Py) = g(x, y)$ for any tangent vectors $x$ and $y$. A. M.
Naveira gave in \cite{Nav} a classification of these manifolds
with respect to the covariant derivative $\n P$, where $\n$  is
the Levi-Civita connection of $g$. This classification is very
similar to the Gray-Hervella classification in \cite{Gray-Her} of
almost Hermitian manifolds. In the paper \cite{Sta-Gri} M.
Staikova and K. Gribachev obtained a classification of the
Riemannian almost product manifolds, for which $\tr P = 0$. In
this case the manifold is even-dimensional. The class $\W_0$ from
the Staikova-Gribachev classification is determined by the
condition $\n P= 0$. A manifold from this class is called a
\emph{Riemannian $P$-manifold}. The class $\W_0$ is an analogue of
the class of K\"ahler manifolds in almost Hermitian geometry.

The geometry of a Riemannian almost product manifold $(M, P, g)$
is a geometry of both structures $g$ and $P$. There are important
in this geometry the linear connections in respect of which the
parallel transport determine an isomorphism of the tangent spaces
with the structures $g$ and $P$. This is valid if and only if
these structures are parallel with respect to such a connection.
In the general case on a Riemannian almost product manifold there
are countless number linear connections regarding which $g$ and
$P$ are parallel. Such connections are called \emph{natural} in
\cite{Mi}.

In the present paper we consider some problems of the geometry of
the natural connections on the class $\W_1$ from the
Staikova-Gribachev classification. This is the class of
\emph{conformal Riemannian $P$-manifolds} or shortly
\emph{$\W_1$-manifolds}. The class $\W_1$ is an analogue of the
class of conformal K\"ahler manifolds in almost Hermitian
geometry.

The paper is organized as follows. In Sec. 2 we give necessary
facts about Riemannian almost product manifolds, the class $\W_1$
and the natural connections on $\W_1$. We recall the notion a
\emph{Riemannian $P$-tensor} on a Riemannian almost product
manifold, which is an analogue of the notion of a K\"ahler tensor
in Hermitian geometry. In Sec. 3 we obtain relations between the
curvature tensors $R$ and $R'$, the Ricci tensors $\rho$ and
$\rho'$ and the scalar curvatures $\tau$ and $\tau'$ of the
Levi-Civita connection $\n$  and a natural connection $\n'$ on a
$\W_1$-manifold. In Sec. 4 we find some necessary and sufficient
conditions the curvature tensor $R'$ of $\n'$ to be a Riemannian
$P$-tensor. The most important result in this section is the
classification \thmref{thm-4.3}. In it there are separated two
special natural connections $D$ and $\widetilde{D}$, whose average
connection is the canonical connection, introduced in \cite{Mi}.
In Sec. 5 we find conditions for a natural connection with
parallel torsion. The important result here is \thmref{thm-5.3}.
It characterizes a $\W_1$-manifold with a natural connection,
which has a parallel torsion and a Riemannian $P$-tensor of
curvature.

\section{Preliminaries}
Let $(M,P,g)$ be a \emph{Riemannian almost product manifold}, \ie
a differentiable manifold $M$ with a tensor field $P$ of type
$(1,1)$ and a Riemannian metric $g$ such that $P^2x=x$,
$g(Px,Py)=g(x,y)$ for any $x$, $y$ of the algebra $\X(M)$ of the
smooth vector fields on $M$. Further $x,y,z,w$ will stand for
arbitrary elements of $\X(M)$ or vectors in the tangent space
$T_pM$ at $p\in M$. The \emph{associated metric} $\widetilde{g}$
of $g$ is determined by $\widetilde{g}(x,y)=g(x,Py)$.

In this work we consider manifolds $(M, P, g)$ with $\tr{P}=0$. In
this case $M$ is an even-dimensional manifold. We assume that
$\dim{M}=2n$.

In \cite{Nav} A.M.~Naveira gives a classification of Riemannian
almost product manifolds with respect to the tensor $F$ of type
(0,3), defined by $ F(x,y,z)=g\left(\left(\nabla_x
P\right)y,z\right), $ where $\n$ is the Levi-Civita connection of
$g$. The tensor $F$ has the properties:
\begin{equation*}
    F(x,y,z)=F(x,z,y)=-F(x,Py,Pz),\qquad
    F(x,y,Pz)=-F(x,Py,z).
\end{equation*}

Using the Naveira classification, in \cite{Sta-Gri} M.~Staikova
and K.~Gribachev give a classification of Riemannian almost
product manifolds $(M,P,g)$ with $\tr P=0$. The basic classes of
this classification are $\W_1$, $\W_2$ and $\W_3$. Their
intersection is the class $\W_0$ of the \emph{Riemannian
$P$-manifolds}, determined by the condition $F=0$ (or equivalently
$\n P=0$) \cite{Sta87}. This class is an analogue of the class of
K\"ahler manifolds in the geometry of almost Hermitian manifolds.

A Riemannian almost product manifold $(M,P,g)$ is a
\emph{Riemannian product manifold} if it has a local product
structure, \ie the structure $P$ is integrable. The Riemannian
product manifolds form the class $\W_1\oplus\W_2$ from the
classification in \cite{Sta-Gri}. This class is an analogue of the
class of Hermitian manifolds in almost Hermitian geometry.

The class $\W_1$ from the Staikova-Gribachev classification
contains the manifolds which are locally conformal equivalent to
Riemannian $P$-manifolds. This class plays a similar role of the
role of the class of the conformal K\"ahler manifolds in almost
Hermitian geometry. We will say that a manifold from the class
$\W_1$ is a \emph{$\W_1$-manifold}.

The characteristic condition for the class $\W_1$ is the following
\begin{equation}\label{2}
\begin{array}{l}
 F(x,y,z)=\frac{1}{2n}\big\{ g(x,y)\ta (z)-g(x,Py)\ta (Pz)
 \big.\\[4pt]
 \phantom{F(x,y,z)=\frac{1}{2n}} +g(x,z)\ta (y)-g(x,Pz)\ta (Py)\big\},
\end{array}
\end{equation}
where the associated 1-form $\ta$ is determined by $
\ta(x)=g^{ij}F(e_i,e_j,x). $ Here $g^{ij}$ will stand for the
components of the inverse matrix of $g$ with respect to a basis
$\{e_i\}$ of $T_pM$ at $p\in M$. The 1-form $\ta$ is
\emph{closed}, \ie $\de\ta=0$, if and only if
$\left(\n_x\ta\right)y=\left(\n_y\ta\right)x$. Moreover, $\ta\circ
P$ is a closed 1-form if and only if
$\left(\n_x\ta\right)Py=\left(\n_y\ta\right)Px$.

In \cite{Sta-Gri} it is proved that
$\W_1=\overline\W_3\oplus\overline\W_6$, where $\overline\W_3$ and
$\overline\W_6$ are the classes from the Naveira classification
determined by the following conditions:
\[
\begin{array}{rl}
\overline\W_3:& \quad F(A,B,\xi)=\frac{1}{n}g(A,B)\ta^v(\xi),\quad
F(\xi,\eta,A)=0,
\\[4pt]
\overline\W_6:& \quad
F(\xi,\eta,A)=\frac{1}{n}g(\xi,\eta)\ta^h(A),\quad F(A,B,\xi)=0,
\end{array}
\]
where $A,B,\xi,\eta\in\X(M)$, $PA=A$, $PB=B$, $P\xi=-\xi$,
$P\eta=-\eta$, $\ta^v(x)=\frac{1}{2}\left(\ta(x)-\ta(Px)\right)$,
$\ta^h(x)=\frac{1}{2}\left(\ta(x)+\ta(Px)\right)$. In the case
when $\tr P=0$, the above conditions for $\overline\W_3$ and
$\overline\W_6$ can be written for any $x,y,z$ in the following
form:
\begin{equation}\label{4'}
\begin{array}{rl}
    \overline\W_3: \quad
    &F(x,y,z)=\frac{1}{2n}\bigl\{\left[g(x,y)+g(x,Py)\right]\ta(z)\\[4pt]
    &+\left[g(x,z)+g(x,Pz)\right]\ta(y)\bigr\},\quad
    \ta(Px)=-\ta(x),
\end{array}
\end{equation}
\begin{equation}\label{5'}
\begin{array}{rl}
    \overline\W_6: \quad
    &F(x,y,z)=\frac{1}{2n}\bigl\{\left[g(x,y)-g(x,Py)\right]\ta(z)\\[4pt]
    &+\left[g(x,z)-g(x,Pz)\right]\ta(y)\bigr\},\quad
    \ta(Px)=\ta(x).
\end{array}
\end{equation}

In \cite{Sta-Gri}, a tensor $L$ of type (0,4) with pro\-per\-ties%
\begin{equation}\label{4}
L(x,y,z,w)=-L(y,x,z,w)=-L(x,y,w,z),
\end{equation}
\begin{equation}\label{5}
L(x,y,z,w)+L(y,z,x,w)+L(z,x,y,w)=0
\end{equation}
is called a \emph{curvature-like tensor}. Such a tensor on a
Riemannian almost product manifold $(M,P,g)$ with the property
\begin{equation}\label{6}
L(x,y,Pz,Pw)=L(x,y,z,w)
\end{equation}
is called a \emph{Riemannian $P$-tensor} in \cite{Mek1}. This
notion is an analogue of the notion of a K\"ahler tensor in
Hermitian geometry.

Let $S$ be a (0,2)-tensor on a Riemannian almost product manifold.
In \cite{Sta-Gri} it is proved that the tensor $\psi_1(S)(x,y,z,w)
=g(y,z)S(x,w)-g(x,z)S(y,w)+S(y,z)g(x,w)-S(x,z)g(y,w)$ is
curvature-like if and only if $S(x,y)=S(y,x)$, and the tensor
$\psi_2(S)(x,y,z,w)=\psi_1(S)(x,y,Pz,Pw)$ is curvature-like if and
only if $S(x,Py)=S(y,Px)$. Obviously
$\psi_2(S)(x,y,Pz,Pw)=\psi_1(S)(x,y,z,w)$. The tensors
$\pi_1=\frac{1}{2}\psi_1(g)$, $\pi_2=\frac{1}{2}\psi_2(g)$,
$\pi_3=\psi_1(\widetilde{g})=\psi_2(\widetilde{g})$ are
curvature-like, and the tensors $\pi_1+\pi_2$, $\pi_3$ are
Riemannian $P$-tensors.

The linear connections in our investigations have a torsion. Let
$\nn$ be a linear connection with a tensor $Q$ of the
transformation $\n \rightarrow\nn$ and a torsion $T$, \ie $\nn_x
y=\n_x y+Q(x,y)$, $T(x,y)=\nn_x y-\nn_y x-[x,y]$. The
corresponding (0,3)-tensors are defined by $Q(x,y,z)=g(Q(x,y),z)$,
$T(x,y,z)=g(T(x,y),z)$. The symmetry of $\n$ implies
$T(x,y)=-T(y,x)=Q(x,y)-Q(y,x)$.

A linear connection $\nn$ on a Riemannian almost product manifold
$(M,P,g)$ is called a \emph{natural connection} if $\nn P=\nn
g=0$. If $\nn$ is a linear connection on $(M,P,g)$, then it is  a
natural connection if and only if \cite{Mi}:
\begin{equation}\label{11'}
    F(x,y,z)=Q(x,y,Pz)-Q(x,Py,z),\quad
    Q(x,y,z)=-Q(x,z,y).
\end{equation}

 The curvature tensor $R$ of $\n$ is determined by
$R(x,y)z=\nabla_x \nabla_y z - \nabla_y \nabla_x z -
    \nabla_{[x,y]}z$ and the corresponding tensor of type (0,4) is defined as
follows $R(x,y,z,w)=g(R(x,y)z,w)$. We denote the Ricci tensor and
the scalar curvature for $\n$ by $\rho$ and $\tau$, respectively,
\ie $\rho(y,z)=g^{ij}R(e_i,y,z,e_j)$ and
$\tau=g^{ij}\rho(e_i,e_j)$. Analogously there are defined the
curvature tensor $R'$ the Ricci tensor $\rho'$ and the scalar
curvature $\tau'$ for any connection $\n'$.

Further $\n'$ will stand for a natural connection on a Riemannian
almost product manifold $(M,P,g)$. Then it is valid the identity
\cite{Dobr1103}
\begin{equation}\label{12}
\begin{split}
R(x,y,z,w)=R'(x,y,z,w)-Q(T(x,y),z,w)- \left(\nn_x Q\right)(y,z,w)\\[4pt]
+\left(\n_yQ\right)(x,z,w)+g\bigl(Q(x,z),Q(y,w)\bigr)-g\bigl(Q(y,z),Q(x,w)\bigr).
\end{split}
\end{equation}

In \cite{Dobr1101} it is established that the natural connections
$\n'$  on a $\W_1$-manifold $(M,P,g)$ form a 2-parametric family,
where the torsion $T$ of $\n'$ on $(M,P,g)$ is determined by
\begin{equation}\label{13}
\begin{array}{l}
    T(x,y,z)=\frac{1}{2n}\left\{g(y,z)\ta(Px)-g(x,z)\ta(Py)\right\}\\[4pt]
            +\lm\left\{g(y,z)\ta(x)-g(x,z)\ta(y)+g(y,Pz)\ta(Px)-g(x,Pz)\ta(Py)\right\}\\[4pt]
            +\mu\left\{g(y,Pz)\ta(x)-g(x,Pz)\ta(y)+g(y,z)\ta(Px)-g(x,z)\ta(Py)\right\}.
\end{array}
\end{equation}

%%%%%%%%%%%%%%%%%%%%%%%%%%%%%%%%%%%%%%%%%%%%%% 3

\section{Curvature properties of natural connections on $\W_1$-manifolds}
According to \cite{Hay}, for the torsion $T$ of $\n'$ and the
tensor $Q$ of the transformation $\n \rightarrow\nn$ it is valid
$2Q(x,y,z)=T(x,y,z)-T(y,z,x)+T(z,x,y)$. Then from \eqref{13} for a
$\W_1$-manifold $(M,P,g)$ holds
\begin{equation}\label{16}
    Q(x,y,z)=T(z,y,x).
\end{equation}
Having in mind \eqref{16}, the equality \eqref{13} implies
\begin{equation}\label{17}
\begin{array}{l}
    Q(y,z,w)=g(y,z)\left\{\lm\ta(w)+\left(\mu+\frac{1}{2n}\right)\ta(Pw)\right\}\\[4pt]
\phantom{Q(y,z,w)}-g(y,w)\left\{\lm\ta(z)+\left(\mu+\frac{1}{2n}\right)\ta(Pz)\right\}\\[4pt]
+g(y,Pz)\left\{\lm\ta(Pw)+\mu\ta(w)\right\}-g(y,Pw)\left\{\lm\ta(Pz)+\mu\ta(z)\right\},
\end{array}
\end{equation}
\begin{equation}\label{18}
\begin{array}{l}
    Q(y,z)=g(y,Pz)\left\{\lm
P\Omega+\mu\Omega\right\}y-\left\{\lm\ta(z)+\left(\mu+\frac{1}{2n}\right)\ta(Pz)\right\}\\[4pt]
\phantom{Q(y,z)}+g(y,z)\left\{\lm\Omega+\left(\mu+\frac{1}{2n}\right)P\Omega\right\}
-\left\{\lm\ta(Pz)+\mu\ta(z)\right\}Py,
\end{array}
\end{equation}
where $\Omega$ is determined by $g(\Omega,x)=\ta(x)$. According to
\eqref{18}, the connection $\n'$, determined by the pair of
parameters $(\lm,\mu)$ has the form
\begin{equation}\label{18'}
\begin{array}{l}
    \n'_xy=\n_xy+g(x,Py)\left\{\lm
P\Omega+\mu\Omega\right\}-\left\{\lm\ta(Py)+\mu\ta(y)\right\}Px\\[4pt]
+g(x,y)\left\{\lm\Omega+\left(\mu+\frac{1}{2n}\right)P\Omega\right\}
-\left\{\lm\ta(y)+\left(\mu+\frac{1}{2n}\right)\ta(Py)\right\}x.
\end{array}
\end{equation}

We denote
   $V(x,w)=\lm\left(\n'_x\ta\right)Pw+\mu\left(\n'_x\ta\right)w$, $U(x,w)=\lm\left(\n'_x\ta\right)w
   +\left(\mu+\frac{1}{2n}\right)\left(\n'_x\ta\right)Pw$.
   Since $\n' P=\n' g=\n'\widetilde{g}=0$, from \eqref{17} it follows
\begin{equation}\label{19}
\begin{array}{l}
   \left(\n'_xQ\right)(y,z,w)=g(y,z)U(x,w)-g(y,w)U(x,z)\\[4pt]
\phantom{\left(\n'_xQ\right)(y,z,w)}+g(y,Pz)V(x,w)-g(y,Pw)V(x,z).
\end{array}
\end{equation}

From \eqref{11'} and \eqref{16} it follows
$Q(T(x,y),z,w)=-g(T(x,y),T(z,w))$ which together with \eqref{18}
and \eqref{19}  we apply in \eqref{12}. By appropriate
calculations and using the  notations
\begin{equation}\label{22}
\begin{array}{l}
    p=\lm\Omega+\left(\mu+\frac{1}{2n}\right)P\Omega,\quad q=\lm
P\Omega+\mu\Omega,
\end{array}
\end{equation}
\begin{equation}\label{23}
\begin{array}{l}
    S'(y,z)=U(y,z)-\frac{1}{2n}\left\{\lm\ta(y)\ta(Pz)+\mu\ta(y)\ta(z)\right\},
\end{array}
\end{equation}
\begin{equation}\label{24}
\begin{array}{l}
    S''(y,z)=V(y,Pz)+\frac{1}{2n}\left\{\lm\ta(Py)\ta(Pz)+\mu\ta(Py)\ta(z)\right\},
\end{array}
\end{equation}
we obtain the following
\begin{thm}\label{thm-3.1}
The following relation is valid
\begin{equation}\label{21}
    R=R'-g(p,p)\pi_1-g(q,q)\pi_2-g(p,q)\pi_3-\psi_1(S')-\psi_2(S'').
\end{equation}
\end{thm}

By the traces of the tensors in \eqref{21} we obtain the following
\begin{cor}\label{cor-3.2}
The following relations are valid
\begin{equation*}
\begin{array}{rl}
    \rho(y,z)&=\rho'(y,z)-\left\{(2n-1)g(p,p)+\tr
    S'-g(q,q)\right\}g(y,z)\\[4pt]
    &-\tr\widetilde{S}g(y,Pz)-2(n-1)S'(y,z)+S''(y,z)+S''(Py,Pz),\\[4pt]
    \tau &=\tau'-2n(2n-1)g(p,p)+2ng(q,q)-2n(2n-1)\tr S'+2\tr S'',
\end{array}
\end{equation*}
where $\widetilde{S}(y,z)=S(y,Pz)$.
\end{cor}

%%%%%%%%%%%%%%%%%%%%%%%%%%%%%%%%%%%  4

\section{Natural connections on $\W_1$-manifolds with Riemannian $P$-tensor of curvature}
Because of $\n'g=\n'P=0$, the properties \eqref{4} and \eqref{6}
are valid for $R'$. Therefore $R'$ is a Riemannian $P$-tensor if
and only if $R'$ satisfies \eqref{5}. This property is valid for
the curvature-like tensors $R$, $\pi_1$, $\pi_2$  and $\pi_3$.
Then, according to \eqref{21}, $R'$ is a Riemannian $P$-tensor if
and only if $\psi_1(S')$ and $\psi_2(S'')$ are curvature-like, \ie
 $S'(y,z)=S'(z,y)$ and $S''(y,Pz)=S''(z,Py)$.
In other words, using \eqref{23} and
 \eqref{24}, it is valid the following
\begin{prop}\label{prop-4.1}
The curvature tensor of $\n'$, determined by \eqref{18'}, is a
Riemannian $P$-tensor if and only if the following conditions are
valid:
\begin{equation}\label{25}
\begin{array}{l}
    U(y,z)-U(z,y)+\frac{\lm}{2n}\left\{\ta(Py)\ta(z)-\ta(y)\ta(Pz)\right\}=0,\\[4pt]
    V(y,z)-V(z,y)+\frac{\mu}{2n}\left\{\ta(Py)\ta(z)-\ta(y)\ta(Pz)\right\}=0.
\end{array}
\end{equation}
\end{prop}

According to \eqref{18}, there are valid the equalities
\begin{equation}\label{26'}
\begin{array}{rl}
    &\left(\n'_y\ta\right)z-\left(\n'_z\ta\right)y=\left(\n_y\ta\right)z-\left(\n_z\ta\right)y\\[4pt]
    &\phantom{\left(\n'_y\ta\right)z-\left(\n'_z\ta\right)y}-\frac{1}{2n}\left\{\ta(Py)\ta(z)-\ta(y)\ta(Pz)\right\},\\[4pt]
        &\left(\n'_y\ta\right)Pz-\left(\n'_z\ta\right)Py=\left(\n_y\ta\right)Pz-\left(\n_z\ta\right)Py
\end{array}
\end{equation}
and then \propref{prop-4.1} implies the following
\begin{prop}\label{prop-4.2}
The curvature tensor of $\n'$, determined by \eqref{18'}, is a
Riemannian $P$-tensor if and only if the following conditions are
valid:
\begin{equation}\label{27}
\begin{array}{l}
    \lm\left\{\left(\n_y\ta\right)z-\left(\n_z\ta\right)y\right\}
    +\left(\mu+\frac{1}{2n}\right)\left\{\left(\n_y\ta\right)Pz-\left(\n_z\ta\right)Py\right\}=0,\\[4pt]
    \mu\left\{\left(\n_y\ta\right)z-\left(\n_z\ta\right)y\right\}
    +\lm\left\{\left(\n_y\ta\right)Pz-\left(\n_z\ta\right)Py\right\}=0.
\end{array}
\end{equation}
\end{prop}

The conditions \eqref{27} form a homogeneous linear system for
$x_1=\left(\n_y\ta\right)z-\left(\n_z\ta\right)y$,
$x_2=\left(\n_y\ta\right)Pz-\left(\n_z\ta\right)Py$ with a
determinant $\Delta=\lm^2-\mu^2-\frac{\mu}{2n}$.

\textbf{Case I}: $\Delta=0$. Then \eqref{27} has a nonzero
solution $(x_1,x_2)$ and for $\lm$ and $\mu$ are possible the
following three subcases: a) $\lm=\mu=0$; b) $\lm=0$,
$\mu=-\frac{1}{2n}$; c) $\lm\neq 0$,
$\lm^2-\mu^2-\frac{\mu}{2n}=0$.

Let us consider the subcase a). If $R'$ is a Riemannian $P$-tensor
then the first equation of \eqref{27} implies that $\ta\circ P$ is
closed. Since \eqref{27} has a nonzero solution then $\ta$ is not
closed. Vice versa, if $\ta\circ P$ is closed, then, according to
\cite{Dobr1103}, $R'$ is a Riemannian $P$-tensor.

Let us consider the subcase b). If $R'$ is a Riemannian $P$-tensor
then the second equation of \eqref{27} implies that $\ta$ is a
closed. Since \eqref{27} has a nonzero solution then $\ta\circ P$
is not closed. Vice versa, if $\ta$ is closed, then the first
equation of \eqref{26'} implies $
\left(\n'_y\ta\right)z-\left(\n'_z\ta\right)y=
\frac{1}{2n}\left\{\ta(y)\ta(Pz)-\ta(Py)\ta(z)\right\} $ and since
$\lm=0$ the second equality of \eqref{27} holds. Obviously, the
first equation of \eqref{27} is valid, too. Then, according to
\propref{prop-4.1}, $R'$ is a Riemannian $P$-tensor.

Let us consider the subcase c). Thus we have also $\mu\neq 0$. If
$R'$ is a Riemannian $P$-tensor, the second equation of \eqref{27}
yields that $\ta$ and $\ta\circ P$ are not closed, since the
system \eqref{27} has a nonzero solution. In this subcase
$\W_1$-manifold $(M,P,g)$ is not belong to the classes
$\overline\W_3$ and $\overline\W_6$. Indeed, in the opposite case
we have $\ta(Pz)=\varepsilon\ta(z)$ $(\varepsilon=\pm 1)$ and
therefore the following equality is valid
\[
    \left(\n_y\ta\right)Pz-\left(\n_z\ta\right)Py
    =\varepsilon\left\{\left(\n_y\ta\right)z-\left(\n_z\ta\right)y\right\}
    -F(y,z,\Omega)+F(z,y,\Omega).
\]
Since \eqref{4'} and \eqref{5'}, it follows
$F(y,z,\Omega)=\frac{1}{2n}\left\{g(y,z)-\varepsilon
g(y,Pz)\right\}\ta(\Omega)$. Using the last two equalities, the
system \eqref{27} gets the form
\begin{equation}\label{27'}
\begin{array}{r}
    \bigl(\lm+\varepsilon\left(\mu+\frac{1}{2n}\right)\bigr)
    \left\{\left(\n_y\ta\right)z-\left(\n_z\ta\right)y\right\}=0,\\[4pt]
    \bigl(\mu+\varepsilon\lm\bigr)
    \left\{\left(\n_y\ta\right)z-\left(\n_z\ta\right)y\right\}=0.
\end{array}
\end{equation}
Because of $\lm\neq 0$, $\lm^2-\mu^2-\frac{\mu}{2n}=0$ in the
considered subcase, we have
$\lm+\varepsilon\left(\mu+\frac{1}{2n}\right)\neq 0$,
$\mu+\varepsilon\lm\neq 0$. Then \eqref{27'} implies that $\ta$ is
closed, which is a contradiction. Therefore
$(M,P,g)\notin\overline\W_3$ and $(M,P,g)\notin\overline\W_6$.

\textbf{Case II}: $\Delta\neq 0$. If $R'$ is a Riemannian
$P$-tensor then the conditions \eqref{27} are satisfied. Since the
homogeneous system has the zero solution $x_1=x_2=0$ in this case
then $\ta$ and $\ta\circ P$ are closed. Vice versa, if $\ta$ and
$\ta\circ P$ are closed, then \eqref{26'} implies the equalities
$\left(\n'_y\ta\right)Pz-\left(\n'_z\ta\right)Py=0$,
$
    \left(\n'_y\ta\right)z-\left(\n'_z\ta\right)y
    =\frac{1}{2n}\left\{\ta(y)\ta(Pz)-\ta(Py)\ta(z)\right\}$
and thus conditions \eqref{25} hold. Then, according to
\propref{prop-4.1}, $R'$ is a Riemannian $P$-tensor.

The obtained results we summarize in the following classification
theorem for the natural connections on a $\W_1$-manifold.
\begin{thm}\label{thm-4.3}
Let $R'$ be the curvature tensor of a natural connection $\n'$
determined by \eqref{18'} on a $\W_1$-manifold $(M,P,g)$. Then the
all possible cases are as follows:
\begin{enumerate}
    \item[i)] If $\n'$ is the connection $D$ determined by
    $\lm=\mu=0$, then $R'$ is a Riemannian $P$-tensor if and only
    if the 1-form $\ta$ is not closed and the 1-form $\ta\circ P$ is
    closed;
    \item[ii)] If $\n'$ is the connection $\widetilde{D}$ determined by
    $\lm=0$, $\mu=-\frac{1}{2n}$, then $R'$ is a Riemannian $P$-tensor if and only
    if the 1-form $\ta$ is closed and the 1-form $\ta\circ P$ is
    not closed;
    \item[iii)] If $\n'$ is a connection for which
    $\lm^2-\mu^2-\frac{\mu}{2n}\neq 0$, then $R'$ is a Riemannian $P$-tensor if and only
    if the 1-forms $\ta$ and $\ta\circ P$ are closed;
    \item[iv)] If $\n'$ is a connection for which $\lm\neq 0$,
    $\lm^2-\mu^2-\frac{\mu}{2n}=0$ and $R'$ is a Riemannian $P$-tensor then the 1-forms $\ta$ and
    $\ta\circ P$ are not closed. In this case $(M,P,g)\notin\overline\W_3$, $(M,P,g)\notin\overline\W_6$.
\end{enumerate}
\end{thm}

Let us remark that the proposition i) of \thmref{thm-4.3} is
proved in \cite{Dobr1103}, where the connection $D$ is
investigated.

In \cite{Mi}, it is introduced a natural connection $\n^c$ on an
arbitrary Riemannian almost product manifold $(M,P,g)$. This
connection, called \emph{canonical}, is an analogue of the
Hermitian connection in Hermitian geometry. Analogues of this
connection are considered on manifolds, corresponding to the class
$\W_1$, in the geometry of almost complex manifolds with Norden
metric (\cite{GanGriMih}) and almost contact manifolds with
B-metric (\cite{ManGri2},\cite{Man4}). In \cite{Dobr1101}, it is
note that $\n^c$ on a Riemannian almost product $\W_1$-manifold is
determined by \eqref{18'} for $\lm=0$, $\mu=-\frac{1}{4n}$.
Therefore, for $\n^c$ it is valid the condition
$\lm^2-\mu^2-\frac{\mu}{2n}\neq 0$ and then the statement iii) of
\thmref{thm-4.3} holds.

Having in mind the values of $\lm$ and $\mu$ determining the
connections $\n^c$, $D$ and $\widetilde{D}$, we obtain the
following
\begin{prop}\label{prop-4.4}
The canonical connection $\n^c$ is the average connection of $D$
and $\widetilde{D}$, \ie $\n^c=\frac{1}{2}(D+\widetilde{D})$.
\end{prop}

%%%%%%%%%%%%%%%%%%%%%%%%%%%%%%%%%%%%%%%%%%%%%%%%%%%%%   5

\section{Natural connection on a $\W_1$-manifold with parallel torsion}
Calculating from \eqref{13} the covariant derivative $\n'T$ and
applying $\n'g=\n'\widetilde{g}=0$, we obtain the following
\begin{prop}\label{prop-5.1}
A natural connection $\n'$ has a parallel torsion if and only if
the 1-form $\ta$ is parallel with respect to $\n'$.
\end{prop}

By \eqref{26'}, \eqref{4'} and \eqref{5'}, from \propref{prop-5.1}
we obtain the following
\begin{cor}\label{cor-5.2}
If a $\W_1$-manifold $(M,P,g)$ admits a natural connection with
parallel torsion then the 1-form $\ta\circ P$ is closed. If
$(M,P,g)$ belongs to the class $\overline\W_3$ or the class
$\overline\W_6$ then the 1-form $\ta$ is closed, too.
\end{cor}

Let $\n'$ be a natural connection with parallel torsion and
Riemannian $P$-tensor of curvature. Then from \propref{prop-5.1}
we have $U=V=0$. In this case \propref{prop-4.1} implies $\lm
W=\mu W=0$, where $W(y,z)=\ta(Py)\ta(z)-\ta(y)\ta(Pz)$. If $W=0$
then we obtain $\ta(Py)=\pm\ta(y)$, \ie $(M,P,g)\in\overline\W_3$
or $(M,P,g)\in\overline\W_6$.
 Thus, we get the following
\begin{thm}\label{thm-5.3}
Let a $\W_1$-manifold $(M,P,g)$ admits a natural connection $\n'$
with parallel torsion and Riemannian $P$-tensor of curvature. Then
the all possible cases are as follows:
\begin{enumerate}
    \item[i)] $\lm=\mu=0$, $W\neq 0$, \ie $\n'=D$,
    $(M,P,g)\notin\overline\W_3$, $(M,P,g)\notin\overline\W_6$;
    \item[ii)] $\lm=\mu=0$, $W= 0$, \ie $\n'=D$,
    $(M,P,g)\in\overline\W_3$ or $(M,P,g)\in\overline\W_6$;
    \item[iii)] $\lm=0$, $\mu\neq 0$, $W= 0$,
    \ie $\n'$ belongs to the 1-parametric family determined by \eqref{18'} for
    $\lm=0$,
    $(M,P,g)\in\overline\W_3$, $(M,P,g)\in\overline\W_6$;
    \item[iv)] $\lm\neq 0$, $\mu=0$, $W= 0$,
    \ie $\n'$ belongs to the 1-parametric family determined by \eqref{18'} for
    $\mu=0$,
    $(M,P,g)\in\overline\W_3$, $(M,P,g)\in\overline\W_6$.
\end{enumerate}
\end{thm}

%%%%%%%%%%%%%%%%%%%%%%%%%%%%%%%%%%%%%%%%%%%%%%%%%%%%%%%%%%%%%%%%%%%%%%%%%%%%

\bigskip
\noindent
\textsl{%D. Gribacheva, D. Mekerov\\
Faculty of Mathematics and Informatics
\\
Paisii Hilendarski University of Plovdiv\\
236 Bulgaria Blvd., 4003 Plovdiv, Bulgaria
\\
e-mails: dobrinka@uni-plovdiv.bg, mircho@uni-plovdiv.bg}

\end{document}